\documentclass[10pt]{article}
 \usepackage{amssymb}
    \usepackage{latexsym}

\title{{\bf Twisted modules for vertex operator algebras and Bernoulli polynomials}}
    \author{B. Doyon\footnote{
B.D. gratefully acknowledges partial support from an NSERC
Postgraduate Scholarship.}, J. Lepowsky \footnote{J.L. and A.M.
gratefully acknowledge partial support from NSF grant
DMS-0070800.} \ and A. Milas \footnotemark[\value{footnote}] }
    \date{}
    
\begin{document}

    \bibliographystyle{alpha}
    \maketitle

    \input amssym.def
    \input amssym

    \newtheorem{rema}{Remark}[section]
    \newtheorem{propo}[rema]{Proposition}
    \newtheorem{theo}[rema]{Theorem}
    \newtheorem{defi}[rema]{Definition}
    \newtheorem{lemma}[rema]{Lemma}
    \newtheorem{corol}[rema]{Corollary}
    \newtheorem{exam}[rema]{Example}

\setcounter{section}{0}
\renewcommand{\theequation}{\thesection.\arabic{equation}}
\renewcommand{\therema}{\thesection.\arabic{rema}}
\setcounter{equation}{0}
\setcounter{rema}{0}

\def\sect#1{\section{#1}\setcounter{equation}{0}\setcounter{rema}{0}}
\def\ssect#1{\subsection{#1}}
\def\sssect#1{\subsubsection{#1}}

\newcommand{\bc}{\begin{center}}
\newcommand{\ec}{\end{center}}
\newcommand{\ba}{\begin{array}}
\newcommand{\ea}{\end{array}}
\newcommand{\nn}{\nonumber \\}
\newcommand{\beq}{\begin{equation}}
\def\sop{x_1^{1/p} \rightarrow \omega_p^s (x_2+x_0)^{1/p}}
\def\rp{x_1^{1/p} \rightarrow \omega_p^r (x_2+x_0)^{1/p}}
\def\xp{x^{1/p} \rightarrow \omega_p^s x^{1/p}}
\newcommand{\eeq}{\end{equation}}
\newcommand{\beqa}{\begin{eqnarray}}
\newcommand{\eeqa}{\end{eqnarray}}
\newcommand{\no}{\nonumber}
\def\bi{\begin{itemize}}
\def\ei{\end{itemize}}
\def\mato#1{\left(\ba{#1}} 
\def\matf{\ea\right)}

\def\eq#1{(\ref{#1})}
\def\lab#1{\label{#1}}

\def\d{\partial}
\def\dv#1#2{\frac{\delta #1}{\delta #2}}
\def\frc#1#2{\frac{#1}{#2}}
\def\b#1{\bar{#1}}
\def\t#1{\tilde{#1}}
\def\h#1{\hat{#1}}
\def\lt#1{\left#1}
\def\rt#1{\right#1}
\def\la{\langle}
\def\ra{\rangle}
\def\<{\langle}
\def\>{\rangle}
\def\:{\mbox{\tiny ${\bullet\atop\bullet}$}}

\def\F{{\Bbb F}}
\def\Z{{\Bbb Z}}
\def\C{{\Bbb C}}
\def\R{{\Bbb R}}
\def\N{{\Bbb N}}
\def\D{{\cal D}}
\def\a#1{{\goth{#1}}}
\def\End{{\rm End\,}}
\def\Hom{{\rm Hom}}
\def\Der{{\rm Der}}
\def\mod{{\rm \,mod\,}}
\def\dim{{\rm dim\,}}
\def\Res{{\rm Res}}
\def\Tr{{\rm Tr}}
\def\proof{{\em Proof: }}
\def\eproof{$\Box$}

\def\for#1{\quad\mbox{#1 }}
\def\com#1{\quad(#1)}
\def\om{\omega_p}

\newcommand{\nordplus}{\mbox{\scriptsize ${+ \atop +}$}}
\newcommand{\nordbullet}{\mbox{\tiny ${\bullet\atop\bullet}$}}

\def\aa{\mbox{\scriptsize ${\ast\atop\ast}$}}
\def\pp{\nordplus}
\def\rr{\nordbullet}
\def\xx{\mbox{\scriptsize ${\times\atop\times}$}}

\def\glim#1{\widehat{\lim_{#1}}}

\def\dpf#1#2{\delta\left(\frac{#1}{#2}\rt)}
\def\pf#1#2{\left(\frac{#1}{#2}\rt)}
\def\id{1_\a{h}}
\sect{Introduction}

This work is a continuation of a series of papers of two of the
present authors \cite{L3}, \cite{L4}, \cite{M1}--\cite{M3}, stimulated
by work of Bloch \cite{Bl}.  In those papers we used the general
theory of vertex operator algebras to study central extensions of
classical Lie algebras and superalgebras of differential operators on
the circle in connection with values of $\zeta$--functions at the
negative integers.  In the present paper, using general principles of
the theory of vertex operator algebras and their twisted modules, we
obtain a bosonic, twisted construction of a certain central extension
of a Lie algebra of differential operators on the circle, for an
arbitrary twisting automorphism.  The construction involves the
Bernoulli polynomials in a fundamental way.  This is explained through
results in the general theory of vertex operator algebras, including a
new identity, which we call ``modified weak associativity.''  This
paper is an announcement. The detailed proofs will appear elsewhere.

More specifically, the present goal is to obtain a new general Jacobi identity for twisted operators, and for related iterates of such operators,
extending the previous analogous results in the untwisted setting in
our papers mentioned above.  As a consequence we obtain twisted
constructions of certain central extensions of Lie algebras of
differential operators on the circle, combining and extending methods
{}from \cite{L3}, \cite{L4}, \cite{M1}--\cite{M3}, \cite{FLM1}, \cite{FLM2} and \cite{DL}.
In those earlier papers we used vertex operator techniques
to analyze untwisted actions of the Lie algebra $\hat{\mathcal{D}}^+$,
studied in \cite{Bl}, on a module for a Heisenberg
Lie algebra of a certain standard type, based on a finite-dimensional
vector space equipped with a nondegenerate symmetric bilinear form.
Now consider an arbitrary isometry $\nu$ of period say $p$, that is,
with $\nu^p={1}$. Here we announce that the corresponding
$\nu$--twisted modules carry an action of the Lie algebra
$\hat{\mathcal{D}}^+$ in terms of twisted vertex
operators, parametrized by certain quadratic vectors in the
untwisted module.

In particular, we extend a result {}from
\cite{FLM1}, \cite{FLM2}, \cite{DL} where actions of the Virasoro algebra were constructed
using twisted vertex operators. In
addition we explicitly compute certain ``correction'' terms for the
generators of the ``Cartan subalgebra'' of $\hat{\mathcal{D}}^+$ that
naturally appear in any twisted construction. These correction terms
are expressed in terms of special values of certain Bernoulli
polynomials. They can in principle be generated, in the theory of
vertex operator algebras, by the formal operator $e^{\Delta_x}$
\cite{FLM1}, \cite{FLM2}, \cite{DL} involved in the construction of a twisted action for a
certain type of
vertex operator algebra. We generate those correction terms in an easier way, using a
new ``modified weak associativity'' relation that is a consequence of
the twisted Jacobi identity.

In \cite{KR} Kac and Radul established a relationship between the Lie algebra
of differential operators on the circle
and the Lie algebra $\widehat{\frak{gl}}(\infty)$; for further work in
this direction, see \cite{AFMO}, \cite{KWY}. Our methods and motivation for studying Lie algebras of differential operators,
based on vertex operator algebras, are new and very different,
so we do not pursue their direction.

\sect{The Lie algebra $\h{\D}^+$ and its untwisted construction}

Let $\mathcal{D}$ be the Lie algebra of formal differential operators on
$\mathbb{C}^\times$ spanned by $t^n D^r$,
where $D=t\,\frac{d}{dt}$ and $n \in \mathbb{Z}$, $r \in \N$ (the nonnegative integers).
Let $\hat{\mathcal{D}}=\mathbb{C}c \oplus \mathcal{D}$ be the
nontrivial one--dimensional central extension (cf. \cite{KR})
with the following commutation relations: \beqa \label{dcom} &&
[t^m f(D), t^n g(D)]= \nn &&
t^{m+n}(f(D+n)g(D)-g(D+m)f(D))+\Psi(t^m f(D),t^ng(D)) c, \no \eeqa
where $f$ and $g$ are polynomials and $\Psi$ is the $2$--cocycle
(cf. \cite{KR}) determined by
$$
    \Psi(t^m f(D), t^n g(D))=-\Psi(t^n g(D), t^m f(D))=
    \delta_{m+n,0} \sum_{i=1}^m f(-i)g(m-i), \; m > 0.
$$
We consider the Lie subalgebra $\D^+$ of $\mathcal{D}$ generated
by the formal differential operators
\beq
\lab{Lnrdef}
    L_{n}^{(r)}=(-1)^{r+1} D^r (t^n D) D^r ,
\eeq
where $n\in\Z,\;r \in\N$ \cite{Bl}.
The subalgebra $\D^+$ has an essentially unique central
extension (cf. \cite{N}) and this extension may be obtained by restriction
of the $2$--cocycle $\Psi$ to $\mathcal{D}^+$.
Let $\hat{\mathcal{D}}^+=\mathbb{C}c \oplus \mathcal{D}^+$
be the nontrivial central extension defined via the slightly
normalized $2$--cocycle $-\frac{1}{2}\Psi$,
and view the elements $L_{n}^{(r)}$ as elements of
$\hat{\mathcal{D}}^+$.
This normalization gives, in particular, the usual Virasoro algebra
bracket relations
$$
    [L_m^{(0)},L_n^{(0)}]=(m-n)L_{m+n}^{(0)}+\frac{m^3-m}{12} \delta_{m+n,0}\,c.
$$

In \cite{Bl} Bloch discovered that the Lie algebra
$\hat{\mathcal{D}}^+$ can be defined in terms of generators that
lead to a simplification of the central term in the Lie bracket relations.
Oddly enough, if we let
\beq\lab{bLnrbloch}
    \bar{L}_n^{(r)}=L_n^{(r)}+ \frac{(-1)^r}{2}\zeta(-1-2r)\delta_{n,0}c,
\eeq
then the central term in the commutator
\beq \label{bl2coc}
    [\bar{L}_{m}^{(r)},\bar{L}_{n}^{(s)}]= \sum_{i={\rm min}(r,s)}^{r+s}
    a_{i}^{(r,s)}(m,n)
    \bar{L}_{m+n}^{(i)}+ \frc{(r+s+1)!^2}{2(2(r+s)+3)!} m^{2(r+s)+3}  \delta_{m+n,0} \, c
\eeq is a pure monomial (here  $a_{i}^{(r,s)}(m,n)$ are structure
constants). In order to conceptualize this simplification
(especially the appearance of $\zeta$--values) one constructs
certain infinite--dimensional projective representations of $\D^+$
using vertex operators.

Let us explain Bloch's construction \cite{Bl}. As in \cite{FLM2},
consider the (infinite--dimensional) Lie algebra $\h{\a{h}}$, the
affinization of an abelian Lie algebra $\a{h}$ of dimension $d$
(over $\mathbb{C}$) with nondegenerate symmetric bilinear form
$\<\cdot,\cdot\>$. The algebra $\h{\a{h}}$ is spanned by
$\alpha(m)\;(\alpha\in\a{h},\;m\in \Z)$ and $C$ (which is
central), satisfying the commutation relations
$$
    [\alpha(m),\beta(n)] = \<\alpha,\beta\> m\delta_{m+n,0} C.
$$
The subalgebra spanned by $\alpha(m)\; (\alpha\in\a{h},\; m\in\Z,\,m\neq 0)$
and $C$ is a Heisenberg Lie algebra, and
$S(\hat{\a{h}}^-)$, where $S(\cdot)$ denotes
the symmetric algebra
and $\hat{\a{h}}^-$ is the span of the $\alpha(m)$ with
$\alpha\in\a{h}$ and $m<0$, carries the structure of an
induced module for $\h{\a{h}}$ with
$\alpha(0)$ acting trivially and $C\mapsto 1$ (cf. \cite{FLM2}).
Then the correspondence
\beq\lab{Lnr}
    L_n^{(r)} \mapsto  \frc12 \sum_{q=1}^d \sum_{j\in\Z} j^r (n-j)^r
\: \alpha_q(j) \alpha_q(n-j) \: \com{n\in\Z}\,,\;
    c\mapsto d,
\eeq
where $\{\alpha_q\}$ is an orthonormal basis of
$\a{h}$ and $\:\cdot\:$ is the usual
normal ordering, which brings $\alpha(n)$ with $n>0$ to the right,
gives a representation of $\hat{\mathcal{D}}^+$.
Let us denote the operator on the right--hand side of (\ref{Lnr}) by
$L^{(r)}(n)$.
In particular, the operators
$L^{(0)}(m) \;(m\in\Z)$ give a well-known
representation of the Virasoro algebra
with central charge $c\mapsto d$,
$$
    [L^{(0)}(m),L^{(0)}(n)] = (m-n)L^{(0)}(m+n) + d \,\frc{m^3-m}{12}\, \delta_{m+n,0},
$$
and the construction \eq{Lnr} for those operators
is the standard realization of the Virasoro algebra on a module for a
Heisenberg Lie algebra (cf. \cite{FLM2}).
The appearance of zeta--values in \eq{bLnrbloch} can be conceptualized
by the following heuristic argument \cite{Bl}:
Suppose that we remove the normal ordering in (\ref{Lnr}) and use
the relation $[\alpha_q(m),\alpha_q(-m)]=m$ to rewrite
$\alpha_q(m) \alpha_q(-m)$, with $m \geq 0$,  as $\alpha_q(-m) \alpha_q(m)+m$.
It is easy to see that the resulting expression
contains an infinite formal divergent series of the form
$$1^{2r+1}+2^{2r+1}+3^{2r+1}+ \cdots .$$
A heuristic argument of Euler's suggests replacing this formal expression
by $\zeta(-1-2r)$, where $\zeta$ is the (analytically continued)
Riemann $\zeta$--function.
The resulting (zeta--regularized) operator is well defined and
gives the action of $\bar{L}_n^{(r)}$; such
operators satisfy the bracket relations (\ref{bl2coc}).

\sect{Twisted modules for vertex operator algebras and the
``modified weak associativity'' relation}

The notion of twisted module for a vertex operator algebra was
formalized in \cite{FFR} and \cite{D} (see also the geometric formulation
in \cite{FrS}), summarizing the basic
properties of the actions of twisted vertex operators discovered in
\cite{FLM1} and \cite{FLM2} (cf. \cite{L1}); the main nontrivial axiom
in this notion is the twisted Jacobi identity of \cite{FLM2} and
\cite{L2} (cf. \cite{FLM1}).  Fix a vertex operator algebra $(V,
Y, {\bf 1}, \omega)$ of central charge $c_V\in\C$ (also called ``rank''),
or $V$ for short (see \cite{FLM2} and \cite{FHL}
for the definition and explicit constructions
of vertex operator algebras and modules, and for necessary ``formal calculus'').
Also fix an automorphism $\nu$ of period $p>0$ of the vertex operator
algebra $V$, that is, a linear automorphism of the vector space $V$
preserving $\omega$ and ${\bf 1}$ such that
\beqa \lab{automnu} &&  \nu Y(v,x)\nu^{-1} = Y(\nu v,x)
\ \mbox{for}\ v\in V, \nn &&
\nu^p=1_V \no
\eeqa
($1_{V}$ being the identity operator on $V$).
In addition, fix a primitive $p$--th root of unity $\om \in \mathbb{C}$.
\begin{defi}\label{tVOA}
{\em A ${\Bbb Q}$-graded $\nu$--twisted $V$--module $M$
is a ${\Bbb Q}$-graded vector space,
$$
    M=\coprod_{n\in {\Bbb Q}}M_{(n)}; \ \mbox{\rm for}\ v\in M_{(n)},\;\mbox{\rm wt}\ v = n,
$$
such that
\beqa
&&  M_{(n)} = 0 \;\; \mbox{ for }\; n \mbox{ sufficiently negative,} \nn
&&  \mbox{\rm dim }M_{(n)}<\infty\;\;\mbox{ for }\; n \in {\Bbb Q} , \no
\eeqa
equipped with a linear map
\begin{eqnarray}
    Y_M(\cdot,x)\,: \ V&\to&(\mbox{\rm End}\; M)[[x^{\frc1p}, x^{-\frc1p}]]\nonumber \\
    v&\mapsto& Y_M(v, x)={\displaystyle \sum_{n\in{\frc1p\Bbb Z}}}v_{n}^\nu x^{-n-1}
    \,,\;\; v_{n}^\nu\in \mbox{\rm End} \; M,
\lab{tvo}
\end{eqnarray}
where $Y_M(v,x)$ is called the {\em twisted vertex operator}
associated with $v$, such that the following conditions hold:\\
{\em truncation condition:}
For every $v
\in V$ and $w \in M$
\beq\lab{truncation}
    v_n^\nu w=0
\eeq
for $n\in\frc1p\Bbb Z$ sufficiently large;\\
{\em vacuum property:}
\begin{equation} \lab{vacuum}
    Y_M({\bf 1}, x)=1_{M};
\end{equation}
{\em Virasoro algebra conditions:}
Let
$$Y_M(\omega,x)=\sum_{n \in \mathbb{Z}} L_M (n)x^{-n-2}.$$
Then
$$[L_M(m),L_M(n)]=(m-n)L_M(m+n)+ c_V\frac{m^3-m}{12}\,\delta_{m+n,0}\, 1_M,$$
$$L_M(0) v=({\rm wt}\ v) v$$ for every homogeneous element $v$, and
$$Y_M(L(-1)u,x)=\frac{d}{dx}Y_M(u,x);$$
{\em Jacobi identity:} For $u,v \in V$,
\beqa\lab{jacobi}
&&  x_0^{-1}\dpf{x_1-x_2}{x_0} Y_M(u,x_1) Y_M(v,x_2)
    - x_0^{-1}\dpf{x_2-x_1}{-x_0} Y_M(v,x_2) Y_M(u,x_1) \no\\
&&  = \frc1p x_2^{-1} \sum_{r=0}^{p-1} \delta\lt(\om^r\lt(\frc{x_1-x_0}{x_2}\rt)^{1/p}\rt)
        Y_M(Y(\nu^r u,x_0)v,x_2),
\eeqa  
where $\delta(x) = \displaystyle{\sum_{n\in{\Bbb Z}} x^n}$
is the {\em formal delta-function}.
}
\end{defi}

Here and below we use the ``binomial expansion convention''---the notational device according to which
binomial expressions are understood to be expanded in nonnegative
integral powers of the second variable.
An important property of a twisted module is that when restricted to the fixed-point subalgebra
$\{v \in V\,|\, \nu v = v\}$, it is a true module: the twisted Jacobi identity \eq{jacobi} reduces to the untwisted
one, as in \cite{FLM2}.

The main commutativity and associativity properties of twisted vertex
operators \cite{L}, together with the fact that these properties are equivalent
to the Jacobi identity, can be reformulated as follows:

\begin{theo} \lab{theostructY}
Let $M$ be a vector space (not assumed to be graded) equipped with a linear
map $Y_M(\cdot,x)$ \eq{tvo} such that the truncation condition \eq{truncation} and the Jacobi
identity \eq{jacobi} hold. Then for $u,\,v\,\in V$ and $w\,\in M$,
there exist $k(u,v) \in\N$ and $l(u,w) \in \frc1p\N$ and
a (non-unique) element $F(u,v,w;x_0,x_1,x_2)$ of
$M((x_0,x_1^{1/p},x_2^{1/p}))$ such that
\beqa
&&  x_0^{k(u,v)} F(u,v,w;x_0,x_1,x_2) \in M[[x_0]]((x_1^{1/p},x_2^{1/p})), \no\\
\lab{propFuvw}
&&  x_1^{l(u,w)} F(u,v,w;x_0,x_1,x_2) \in M[[x_1^{1/p}]]((x_0,x_2^{1/p}))
\eeqa
and
\beqa
    Y_M(u,x_1) Y_M(v,x_2)w &=& F(u,v,w;x_1-x_2,x_1,x_2), \no\\
    Y_M(v,x_2) Y_M(u,x_1)w &=& F(u,v,w;-x_2+x_1,x_1,x_2), \no\\
\lab{structY}
    Y_M(Y(\nu^{-s}u,x_0)v,x_2)w &=&
\lim_{x_1^{1/p} \rightarrow \omega_p^s (x_2+x_0)^{1/p}}
 F(u,v,w;x_0,x_1,x_2)
\eeqa
for $s\in\Z$ (where we are using the binomial expansion convention).
Conversely, let $M$ be a vector space equipped with a linear map $Y_M(\cdot,x)$ \eq{tvo} such that the truncation
condition \eq{truncation} and the preceding statement hold, except that $k(u,v)$
($\in\N$) and $l(u,w)$ ($\in \frc1p\N$) may depend on all three of $u,v$ and $w$.
Then the Jacobi identity \eq{jacobi} holds.
\end{theo}

It is important to note that since $k(u,v)$ can be greater than 0,
the formal series $F(u,v,w;x_1-x_2,x_1,x_2)$ and
$F(u,v,w;-x_2+x_1,x_1,x_2)$ are not in general equal. The formal
limit procedure $\displaystyle{\lim_{x_1^{1/p} \rightarrow
\omega_p^s (x_2+x_0)^{1/p}}F(u,v,w;x_0,x_1,x_2)}$ above means that
one replaces each integral power of the formal variable
$x_1^{1/p}$ in the formal series $F(u,v,w;x_0,x_1,x_2)$ by the
corresponding power of the formal series $\omega_p^s
(x_2+x_0)^{1/p}$ (defined using the binomial expansion
convention).

Along with (\ref{propFuvw}), the first two equations of (\ref{structY})
represent what we call ``formal commutativity'' for twisted vertex operators,
while the first and last equations of (\ref{structY}) represent ``formal associativity'' for twisted
vertex operators. When specialized to the untwisted case $p=1$ ($\nu=1_V$), these two relations
lead respectively to the usual ``formal commutativity'' and ``formal
associativity'' for vertex operators, as formulated in
\cite{LL} (see also \cite{FLM2} and \cite{FHL}). The first and last
relations of the theorem also imply a significant new relation,
interesting even for the untwisted case,
which we call ``modified weak associativity'':
\begin{theo}\lab{theomwa}
With $M$ and $k(u,v)$ as in Theorem \ref{theostructY},
\beqa
&&  \lim_{x_1^{1/p} \to \om^{s}(x_2+x_0)^{1/p}} \lt( (x_1-x_2)^{k(u,v)} Y_M(u,x_1) Y_M(v,x_2) \rt) \no\\
\lab{modwasst}
&& \qquad=  x_0^{k(u,v)}Y_M(Y(\nu^{-s}u,x_0)v,x_2)
\eeqa
for $u,\,v\,\in V$ and $s\in \Z$.
\end{theo}

\begin{rema}
{\em
The specialization of Theorems \ref{theostructY} and \ref{theomwa} to the case $p=1$ ($\nu=1_V$) and $M=V$
gives statements describing structures and relations for the vertex operator map $Y(\cdot,x)$ in the
vertex (operator) algebra $V$.
}
\end{rema}

\sect{Homogeneous twisted vertex operators and associated relations}

Homogeneous vertex operators are defined by
\beq\lab{homo}
    X(v,x) = Y(x^{L(0)}v,x)
\eeq
for $v \in V$ (see \cite{FLM2}).
{}From the Jacobi identity for vertex operators, it is not hard to obtain
a Jacobi identity for homogeneous vertex
operators \cite{L3}, \cite{L4}, \cite{M1}:
\beqa
    && x_0^{-1}\delta\lt(\frc{x_1-x_2}{x_0}\rt) X(u,x_1)X(v,x_2) -
    x_0^{-1}\delta\lt(\frc{x_2-x_1}{-x_0}\rt) X(v,x_2)X(u,x_1) \no\\
\lab{jacobiX}
    && \  = x_1^{-1} \delta\lt( e^{y} \frc{x_2}{x_1}\rt)
        X(Y[u,y] v,x_2)
\eeqa
where
$$
    y = \log\lt(1+\frc{x_0}{x_2}\rt).
$$
The expressions $\log(1+x)$ and $e^y$, where $x$ and $y$ are formal variables, are
defined by their series expansions in nonnegative powers of $x$ and $y$, respectively.
The use of the formal variable $y$ is natural here in particular because of the appearance
of the vertex operator $Y[u,y]$ defined and studied in \cite{Z1}, \cite{Z2}:
\beq\lab{Zhu}
    Y[u,y] = Y(e^{yL(0)}u,e^y-1).
\eeq
These operators give a new vertex operator algebra isomorphic to $V$,
and the isomorphism corresponds geometrically to a
change-of-coordinates transformation expressed formally as $y\mapsto
e^y-1$ (\cite{Z1}, \cite{Z2}; see \cite{H1}, \cite{H2} for the
generalization to arbitrary coordinate changes).
The Jacobi
identity \eq{jacobiX} thus suggests that there is a close relationship between
homogeneous vertex operators and ``cylindrical coordinates''; we will
comment more on this relationship in an extended version of this announcement.

{}From the Jacobi identity \eq{jacobiX}, one can obtain the commutator
formula (\cite{L3}, \cite{L4}):
\beq \lab{comformnt}
    [X(u,x_1),X(v,x_2)] = \Res_y \delta\lt( e^y\frc{x_2}{x_1}\rt)
        X(Y[ u,y] v,x_2),
\eeq
using a general fact concerning formal series:
$$
    \Res_x h(x) = \Res_y \lt( h(F(y))\,\frc{d}{dy} F(y)\rt)  \for{for} h(x)\in A((x)),\;
    F(x)\in xA[[x]]
$$
where $A$ is a commutative associative
algebra (or more generally, a module for it) and where the
coefficient of $x^1$ in $F(x)$ is invertible.

{}From now on we fix a ${\Bbb Q}$-graded $\nu$-twisted
$V$-module $M$.
We define homogeneous twisted vertex operators by a simple twisted generalization of \eq{homo}:
$$
    X_M(u,x) = Y_M(x^{L(0)}u,x),
$$
as in \cite{FLM2}.
We now state the following twisted generalizations of the results
above: a Jacobi identity, a commutator formula
similar to \eq{comformnt}, and formal commutativity and associativity
properties, including ``modified weak associativity,'' for these
homogeneous twisted vertex operators.

Using the twisted Jacobi identity \eq{jacobi} and the definitions
\eq{homo} and \eq{Zhu}, one can obtain the
twisted generalization of \eq{jacobiX}:
\begin{theo}
For $u,v \in V$,
\beqa \label{commtwi}
    && x_0^{-1}\delta\lt(\frc{x_1-x_2}{x_0}\rt) X_M(u,x_1)X_M(v,x_2) -
    x_0^{-1}\delta\lt(\frc{x_2-x_1}{-x_0}\rt) X_M(v,x_2)X_M(u,x_1) \no\\
\lab{jacobiXt}
    && \  = \frc1p x_1^{-1} \sum_{r=0}^{p-1} \delta\lt( \om^{-r} \lt(e^y \frc{x_2}{x_1}\rt)^{1/p}\rt)
        X_M(Y[\nu^r u,y] v,x_2)
\eeqa
where
\beq
    y = \log\lt(1+\frc{x_0}{x_2}\rt).
\eeq
\end{theo}
Using techniques similar to those in the untwisted case, this Jacobi
identity leads to the commutator formula:
\begin{corol}
For $u,v \in V$,
\beq \lab{comform}
    [X_M(u,x_1),X_M(v,x_2)] = \Res_y \frc1p \sum_{r=0}^{p-1} \delta\lt( \om^{-r} \lt(e^y\frc{x_2}{x_1}\rt)^{1/p} \rt)
        X_M(Y[\nu^r u,y] v,x_2).
\eeq
\end{corol}

Applying $\displaystyle{\lim_{x_0 \to (e^y-1)x_2}}$ to both sides
of the ``modified weak associativity'' relation \eq{modwasst} yields
``modified weak associativity'' for homogeneous twisted vertex
operators:
\begin{theo}
For $u,v \in V$, $s \in \mathbb{Z}$, and
$k(u,v)$ as in Theorem \ref{theostructY}, 
\beqa \lab{modwassht}
&&  \lim_{x_1^{1/p} \to \omega_p^s (e^y x_2)^{1/p}}
\lt(\lt(\frc{x_1}{x_2}-1\rt)^{k(u,v)} X_M(u,x_1)X_M(v,x_2)\rt) \nn
&&\qquad = \lt(e^y-1\rt)^{k(u,v)} X_M(Y[\nu^{-s}u,y]v,x_2).
\eeqa
\end{theo}

Finally, the two first equations of \eq{structY} can be rewritten in terms
of homogeneous vertex operators, and modified weak
associativity \eq{modwassht} can be used to obtain the homogeneous counterpart of the third equation
of \eq{structY}:
\begin{theo}
For $u,\,v\,\in V$, $w\in M$ and $s\in\Z$, we have
\beqa
    X_M(u,x_1) X_M(v,x_2)w &=& G(u,v,w;x_1-x_2,x_1,x_2), \no\\
    X_M(v,x_2) X_M(u,x_1)w &=& G(u,v,w;-x_2+x_1,x_1,x_2), \no\\
\lab{structX}
    X_M(Y[\nu^{-s}u,y]v,x_2)w &=&
    \lim_{ x_1^{1/p} \to \omega_p^s (e^y x_2)^{1/p}}
    G(u,v,w;(e^{y}-1)x_2,x_1,x_2),
\eeqa
where
$$
    G(u,v,w;x_0,x_1,x_2) \in M((x_0,x_1^{1/p},x_2^{1/p}))
$$
and
\beqa
&&  x_0^{k(u,v)} G(u,v,w;x_0,x_1,x_2) \in M[[x_0]]((x_1^{1/p},x_2^{1/p})), \no\\
\lab{propGuvw}
&&  x_1^{l'(u,w)} G(u,v,w;x_0,x_1,x_2) \in M[[x_1^{1/p}]]((x_0,x_2^{1/p})),
\eeqa
for some $k(u,v) \in\N$ and $l'(u,w) \in \frc1p\N$. Here $k(u,v)$
can be taken to be the same as in Theorem \ref{theostructY}.
\end{theo}

Again, along with (\ref{propGuvw}), the first two equations of
(\ref{structX}) represent what we call ``formal commutativity'' for
homogeneous twisted vertex operators, while the first and last
equations of (\ref{structX}) represent ``formal associativity'' for
homogeneous twisted vertex operators. These two formal relations are
also of interest in the untwisted case $p=1$ ($\nu=1_V$), and of
course we have the same relations for the vertex operator algebra $V$.

\sect{Commutator formula for iterates on twisted modules}

``Modified weak associativity'' for homogeneous twisted vertex
operators turns out to be very useful.  Formal limit operations
respect products, under suitable conditions, and using this principle,
one can compute, for instance, commutators of certain iterates in a
natural way.  For our applications, an important commutator is
$$
    [X_M(Y[u_1,y_1]v_1,x_1),X_M(Y[u_2,y_2]v_2,x_2)],
$$
which we would like to express in terms of similar iterates.
Using ``modified weak associativity'' \eq{modwassht} and
the commutator formula \eq{comform} one can prove:
\begin{theo}
For $u_1,v_1,u_2,v_2 \in V$,
\label{lastf}
\beqa \no
&&  [X_M(Y[u_1,y_1]v_1,x_1),X_M(Y[u_2,y_2]v_2,x_2)] = \\
\no&&   \Res_y \frc1p \sum_{r=0}^{p-1}
\\\no&& \lt\{
        \delta\lt(\om^{r}\lt(e^{-y_2-y}\frc{x_1}{x_2}\rt)^{1/p}\rt)
        X_M(Y[u_1,y_1+y]Y[\nu^{-r} v_2,-y_2]Y[v_1,y]\nu^{-r}u_2,e^{-y}x_1)
    \rt.
\\\no&& +
        \delta\lt(\om^{r}\lt(e^{-y}\frc{x_1}{x_2}\rt)^{1/p}\rt)
        X_M(Y[u_1,y_1+y]Y[\nu^{-r}u_2,y_2] Y[v_1,y]\nu^{-r} v_2,e^{-y}x_1)
\\\no&& +
        \delta\lt(\om^{r}\lt(e^{y_1-y_2-y}\frc{x_1}{x_2}\rt)^{1/p}\rt)
        X_M(Y[Y[\nu^{-r}v_2,-y_2]Y[u_1,y]\nu^{-r}u_2,y_1-y]v_1,x_1)
\\\no&& \lt.    +
        \delta\lt(\om^{r}\lt(e^{y_1-y}\frc{x_1}{x_2}\rt)^{1/p}\rt)
        X_M(Y[Y[\nu^{-r}u_2,y_2]Y[u_1,y]\nu^{-r}v_2,y_1-y]v_1,x_1)
    \rt\}.
\eeqa
\end{theo}

This generalizes the main commutator formula in \cite{L4} and is
related to a similar commutator formula in \cite{M1}--\cite{M3}.

\sect{Main results}

We now obtain a representation of $\hat{\mathcal{D}}^+$ on a certain
natural module for a twisted affine Lie algebra based on a
finite--dimensional abelian Lie algebra (essentially a twisted
Heisenberg Lie algebra), generalizing Bloch's untwisted representation
on the module $S(\hat{\a{h}}^-)$ constructed in Section 2.  This is
also a generalization of the twisted Virasoro algebra construction
(see \cite{FLM1}, \cite{FLM2}, \cite{DL}).

Let $\a{h}$ be a finite--dimensional abelian Lie algebra (over $\C$)
on which there is a nondegenerate symmetric bilinear form
$\la\cdot,\cdot\ra$.  Let $\nu$ be an isometry of $\a{h}$ of period
$p>0$:
$$
    \la\nu\alpha,\nu\beta\ra = \la\alpha,\beta\ra ,\qquad
    \nu^p\alpha = \alpha
$$
for all $\alpha,\beta\in\a{h}$.  We assume that $\nu$ preserves a
rational lattice in $\a{h}$.  One knows that $S(\hat{\a{h}}^-)$
carries a natural structure of vertex operator algebra, with
${\bf 1}=1$, and that $\nu$
lifts naturally to an automorphism, which we continue to call $\nu$,
of period $p$ of $S(\hat{\a{h}}^-)$ (cf. \cite{FLM2}).  We proceed as
in \cite{L1}, \cite{FLM1}, \cite{FLM2} and \cite{DL} to construct a
space $S[\nu]$ that carries a natural structure of $\nu$--twisted
module for $S(\hat{\a{h}}^-)$.

Recalling our primitive $p$--th root of unity $\omega_p$, for $r \in
\Z$ set
$$
\a{h}_{(r)} = \{\alpha\in\a{h} \;|\; \nu\alpha = \om^r\alpha\}
    \subset \a{h}.
$$
For $\alpha\in\a{h}$, denote by $\alpha_{(r)},\; r\in\Z$, its
projection on $\a{h}_{(r)}$.
Define the $\nu$-twisted affine Lie algebra $\hat{\a{h}}[\nu]$ associated with the abelian Lie algebra $\a{h}$ by
\beq
    \hat{\a{h}}[\nu] = \coprod_{n\in\frc1p\Z} \a{h}_{(pn)} \otimes t^n \oplus \C C
\eeq
with
\beqa
&&  [\alpha\otimes t^m,\beta\otimes t^n] = \<\alpha,\beta\>
 m\delta_{m+n,0}\, C
 \com{\alpha\in\a{h}_{(pn)},\; \beta\in\a{h}_{(pm)},\; m,n\in\frc1p\Z} \no\\
&&  [C,\hat{\a{h}}[\nu]] = 0 .
\eeqa
Set
\beq
    \hat{\a{h}}[\nu]^+ = \coprod_{n>0} \a{h}_{(pn)}\otimes t^n ,\qquad \hat{\a{h}}[\nu]^- = \coprod_{n<0} \a{h}_{(pn)} \otimes t^n.
\eeq
The subalgebra
\beq
    \hat{\a{h}}[\nu]_{\frc1p\Z} = \hat{\a{h}}[\nu]^+ \oplus \hat{\a{h}}[\nu]^- \oplus \C C
\eeq
is a Heisenberg Lie algebra. Form the induced (level-one) $\hat{\a{h}}[\nu]$-module
\beq
    S[\nu] = \mathcal{U}(\hat{\a{h}}[\nu]) \otimes_{\mathcal{U}\lt( \hat{\a{h}}[\nu]^+ \oplus \a{h}_{(0)} \oplus \C C\rt)} \C
        \simeq S(\hat{\a{h}}[\nu]^-) \for{(linearly),}
\eeq
where $\hat{\a{h}}[\nu]^+\oplus \a{h}_{(0)}$ acts trivially on $\C$
and $C$ acts as 1; $\mathcal{U}(\cdot)$ denotes universal enveloping
algebra.  Then $S[\nu]$ is irreducible under
$\hat{\a{h}}[\nu]_{\frc1p\Z}$.
We will use the notation $\alpha(n)\;(\alpha\in\a{h}_{(pn)},\,
n\in\frc1p\Z)$ for the action of
$\alpha\otimes t^n \in \hat{\a{h}}[\nu]$ on $S[\nu]$.
As we mentioned above, the $\hat{\a{h}}[\nu]$-module $S[\nu]$ is
naturally a $\nu$--twisted module for the vertex operator algebra
$S(\hat{\a{h}}^-)$, and its structure and general properties are
important in establishing our results described below.

\begin{rema}
{\em The special case where $p=1$ ($\nu=1_{\a{h}}$) corresponds to the
$\hat{\a{h}}$-module $S(\hat{\a{h}}^-)$ discussed in Section 2.
}
\end{rema}

We now consider certain operators on $S[\nu]$ that we will use for constructing a representation
of $\hat{\mathcal{D}}^+$. For $\alpha \in \goth{h}$ we define the following homogeneous twisted vertex operator
acting on $S[\nu]$:
$$
    \alpha^{\nu}\langle x\rangle=X_{S[\nu]}(\alpha(-1){\bf 1},x)=\sum_{n \in \frac{1}{p}\mathbb{Z}} \alpha(n)x^{-n}.
$$
Choosing an orthonormal basis $\{\alpha_q|q=1,\ldots,d\}$ of $\goth{h}$, we define
the following two formal series acting on $S[\nu]$:
\beqa \lab{Lnuy1y2x}
    L^{\nu;y_1,y_2}\<x\>&=& \frc12 \sum_{q=1}^{d}
\: \alpha^\nu_q\langle e^{y_1}x\rangle \alpha^\nu_q\langle e^{y_2}x\rangle  \: \no \\
&& - \frc12 \frc{\d}{\d y_1} \lt(\sum_{k=0}^{p-1} \frac{e^{\frac{k(-y_1+y_2)}{p}}
{\rm dim} \ \goth{h}_{(k)} - 1}
{1-e^{-y_1+y_2}} \rt)
\eeqa
and
\beqa \lab{bLnuy1y2x}
    {\bar{L}}^{\nu;y_1,y_2}\<x\>&=& \frc12 \sum_{q=1}^{d} \: \alpha^\nu_q\langle e^{y_1}x\rangle
\alpha^\nu_q\langle e^{y_2}x\rangle  \: \no \\
&& - \frc12 \frc{\d}{\d y_1} \lt(\sum_{k=0}^{p-1} \frac{e^{\frac{k(-y_1+y_2)}{p}}
{\rm dim} \ \goth{h}_{(k)}}
{1-e^{-y_1+y_2}} \rt).
\eeqa
\begin{rema}
{\em
In the special case $p=1$ and $d=1$, the operators
$L^{\nu;y_1,y_2}\<x\>$ and ${\bar L}^{\nu;y_1,y_2}\<x\>$, respectively, specialize to the
operators $L^{(y_1,y_2)}(x)$ and ${\bar L}^{(y_1,y_2)}(x)$ of \cite{L3}, \cite{L4}.
}
\end{rema}

The formal series \eq{bLnuy1y2x} can be rewritten in the following form:
$$
    \b{L}^{\nu;y_1,y_2}\<x_2\>
    =
    \frc12 \lim_{x_1 \to  x_2} \sum_{q=1}^d \lt(\lt( \frac{ \frac{x_1}{x_2}e^{y_1-y_2} -1 }
{e^{y_1-y_2} -1} \rt)^k \alpha^\nu_q \langle e^{y_1}x_1\rangle \alpha^\nu_q\langle e^{y_2}x_2\rangle \rt)
$$
for any fixed $k\in\N,\,k\ge2$.
Using ``modified weak associativity'' \eq{modwassht}, we can then identify this with a particular
iterate of vertex operators:
\begin{propo} \label{mathcliter}
With the notation as above,

\beq \label{iterbar}
    {\bar{L}}^{\nu;y_1,y_2}\<x\>= X_{S[\nu]} \lt(\frac{1}{2} \sum_{q=1}^d Y[\alpha_q(-1) {\bf 1} ,y_1-y_2]
    \alpha_q(-1) {\bf 1} ,e^{y_2}x \rt).
\eeq
\end{propo}

Once this identification is made, our twisted construction of $\hat{\mathcal{D}}^+$ is a simple consequence
of the general theory of twisted modules for vertex operator algebras. In particular, consider
the commutator formula for the untwisted operators ${\bar L}^{(y_1,y_2)}(x)$ announced in \cite{L3}
(for a proof see \cite{M1}). Along with Proposition \ref{mathcliter} and the general properties of twisted modules
for vertex operator algebras, it implies the same commutator formula for the twisted operators
${\bar L}^{\nu;y_1,y_2}\<x\>$:
\begin{propo}
With the notation as above,
\begin{eqnarray}\label{Lbarbracketsalg}
\lefteqn{[{\bar L}^{\nu;y_1,y_2}\<x_1\>,{\bar
L}^{\nu;y_3,y_4}\<x_2\>]=}\nonumber \\
&&= - {\frac{1}{2}} \frac{\partial}{\partial y_1} \biggl({\bar
L}^{\nu;-y_1+y_2+y_3,y_4}\<x_2\> \delta
\left({\frac{e^{y_1}x_1}{e^{y_3}x_2}}\right)+ {\bar
L}^{\nu;-y_1+y_2+y_4,y_3}\<x_2\>
\delta \left({\frac{e^{y_1}x_1}{e^{y_4}x_2}}\right)\biggr)\nonumber\\
&&- {\frac{1}{2}} \frac{\partial}{\partial y_2} \biggl({\bar
L}^{\nu;y_1-y_2+y_3,y_4}\<x_2\> \delta
\left({\frac{e^{y_2}x_1}{e^{y_3}x_2}}\right)+ {\bar
L}^{\nu;y_1-y_2+y_4,y_3}\<x_2\> \delta
\left({\frac{e^{y_2}x_1}{e^{y_4}x_2}}\right)\biggr).
\end{eqnarray}
\end{propo}

In fact, this commutator formula can alternatively be directly obtained from our general commutator formula,
Theorem \ref{lastf}.
Formula \eq{Lbarbracketsalg} provides a representation of the Lie algebra $\hat{\mathcal{D}}^+$
on the $\hat{\a{h}}[\nu]$-module $S[\nu]$ via our twisted operators, generalizing the untwisted case. More precisely, let
\beqa
    {L}^{\nu;y_1,y_2}\<x\>&=&
        \sum_{n\in\Z,\;r_1,r_2\,\in\N} {L}^{(r_1,r_2)}(n) x^{-n} \frc{y_1^{r_1}y_2^{r_2}}{r_1!r_2!}, \no \\
    \b{L}^{\nu;y_1,y_2}\<x\>&=& \frc{1}2 \frc{d}{(y_1-y_2)^2} +
        \sum_{n\in\Z,\;r_1,r_2\,\in\N} \b{L}^{(r_1,r_2)}(n) x^{-n}
        \frc{y_1^{r_1}y_2^{r_2}}{r_1!r_2!}. \nonumber
\eeqa
Then the following holds (recall the generators \eq{Lnrdef}, \eq{bLnrbloch} of $\hat{\mathcal{D}}^+$):
\begin{theo} \label{main1}
Let
\beqa
    {L}^{(r)}(n) &=& {L}^{(r,r)}(n) \com{n\in\Z,\,r\in\N}, \no \\
    \b{L}^{(r)}(n) &=& \b{L}^{(r,r)}(n) \com{n\in\Z,\,r\in\N}. \no
\eeqa
\begin{itemize}
\item[(a)] The assignment
$$ L_n^{(r)} \mapsto L^{(r)}(n), \ \ c \mapsto d,$$
defines a representation of the Lie algebra
$\hat{\mathcal{D}}^+$ on $S[\nu]$.
\item[(b)]
The assignment
$$\bar{L}_n^{(r)} \mapsto \bar{L}^{(r)}(n), \ \ c \mapsto d $$
also defines a representation of the Lie algebra $\hat{\mathcal{D}}^+$, with
the central term being a pure monomial, as in \eq{bl2coc}.
\end{itemize}
\end{theo}

Explicit expressions for the operators $L^{(r)}(n)$ and $\b{L}^{(r)}(n)$,
involving Bernoulli polynomials, are easy to obtain from \eq{Lnuy1y2x} and \eq{bLnuy1y2x}:
\beqa \label{lnr}
L^{(r)}(n)&=& \frc12 \sum_{q=1}^{d} \sum_{j\in \frac{1}{p}\Z} j^{r}(n-j)^{r}\,\: \alpha_q
(j)\alpha_q(n-j) \:
\no\\   &&  - \delta_{n,0}\, \frc{(-1)^r}{4(r+1)} \sum_{k=0}^{p-1} \dim\a{h}_{(k)}
        \lt( B_{2(r+1)}(k/p) - B_{2(r+1)} \rt)
\eeqa
and
\beqa \label{opLn}
    \b{L}^{(r)}(n) &=& \frc12 \sum_{q=1}^d \sum_{j\in \frac{1}{p} \Z} j^{r}(n-j)^{r}\,
\: \alpha_q (j)\alpha_q (n-j) \:
\no\\   &&  - \delta_{n,0} \, \frc{(-1)^r}{4(r+1)} \sum_{k=0}^{p-1} \dim\a{h}_{(k)} B_{2(r+1)}(k/p).
\eeqa
{}From our construction, the appearance of Bernoulli polynomials is seen to be directly
related to general properties of homogeneous twisted vertex operators.

The next result is a simple consequence of Theorem \ref{main1}. It describes
the action of the ``Cartan subalgebra'' of $\hat{\mathcal{D}}^+$
on a highest weight vector of a canonical quasi-finite
$\hat{\mathcal{D}}^+$--module; here we are using the terminology of
\cite{KR}.
This corollary gives the ``correction'' terms referred to in the introduction.
\begin{corol}
Given a highest weight $\hat{\mathcal{D}}^+$--module $W$, let
$\delta$ be the linear functional on the ``Cartan subalgebra''
of $\hat{\mathcal{D}}^+$ (spanned by $L^{(k)}_{0}$ for $k\in \N$) defined by
$$
    {L}^{(k)}_{0} \cdot w =(-1)^k\delta\lt( {L}^{(k)}_{0} \rt)w,
$$
where $w$ is a generating highest weight vector of $W$, and let $\Delta(x)$
be the generating function
$$\Delta(x)=\sum_{k \ge 1} \frac{\delta({L}^{(k)}(0))x^{2k}}{(2k)!}$$
(cf. \cite{KR}).
Then for every automorphism $\nu$ of period $p$ as above,
$$\mathcal{U}(\hat{\mathcal{D}}^+) \cdot 1 \subset S[\nu]$$
is a quasi--finite highest weight $\hat{\mathcal{D}}^+$--module satisfying
$$\Delta(x) = \frac{1}{2} \frac{d}{dx} \sum_{k=0}^{p-1} \frac{ e^{\frac{kx}{p}}
{\rm dim} \ \goth{h}_{(k)} -1} {1-e^x}.$$
\end{corol}

Finally, we have an additional result (for the untwisted bosonic case an equivalent result was
obtained in \cite{Bl} and for the spinor constructions in \cite{M2}):
\begin{corol}
The generating function
\beq \label{generators}
X_{S[\nu]}(\sum_{q=1}^d \alpha_q(-m-1)\alpha_q(-m-1){\bf 1},x),
\eeq
$m \in \mathbb{N}$,
defines the same $\hat{\mathcal{D}}^+$--module
as in Theorem \ref{main1}.  That is, every operator
$L^{(r)}(n)$ (or equivalently $\bar{L}^{(r)}(n)$) can be expressed as
a linear combination of the expansion coefficients of the operator
(\ref{generators}), and vice versa.
\end{corol}

\small{

}

\noindent
{\small \sc Department of Physics, Rutgers University, Piscataway,
NJ 08854} \\
{\em E-mail address}: doyon@physics.rutgers.edu, \\
{\small \sc Department of Mathematics, Rutgers University, Piscataway,
NJ 08854} \\
{\em E-mail address}: lepowsky@math.rutgers.edu \\
{\small \sc Department of Mathematics, University of Arizona, Tucson,
AZ 85721} \\
{\em E-mail address}: milas@math.arizona.edu

\end{document}